\documentclass[11pt]{amsart}
\usepackage{amscd,amssymb,palatino}
\usepackage[utf8]{inputenc}
\usepackage[english]{babel}
\usepackage{caption}
\usepackage{etex}
\usepackage[leqno]{amsmath}
\usepackage{amssymb}
\usepackage{color}
\usepackage{shadow}
\usepackage{epsfig}
\usepackage{epic}
\usepackage{eepic}
\usepackage{graphics}
\usepackage{graphicx}
\usepackage{psfrag}
\usepackage{calc}
\usepackage{tikz}
\usepackage[dvipsnames,prologue,table]{pstricks}
\usepackage{pstcol}
\usetikzlibrary{arrows,decorations,patterns,positioning,automata,shadows,fit,shapes,calc,decorations.markings,backgrounds,scopes,decorations.text}

\usepackage{tipa}

\usepackage{fancyhdr}
\fancypagestyle{plain}{\fancyhf{}
\fancyfoot[C]{{\large\sf\bfseries\thepage}}}
\usepackage{calc}
\newcommand\1{\lower 9pt\hbox{\underbar{}}}
\numberwithin{equation}{section}

\newtheorem {Theorem}                   {Theorem}

\newtheorem {RefTheorem}[equation]      {Theorem}

\newtheorem {Conjecture}[equation]      {Conjecture}

\theoremstyle{definition}
\newtheorem {Definition}[equation]{Definition}
\newtheorem {Remark}[equation]          {Remark}

\def\cL{\mathcal L}
\newcommand{\pr} {\smallskip\noindent{\bf Proof\,\,}}
\newenvironment{Proof}  {\pr}{\hspace*{\fill}\qed\\}

\begin{document}

\title[Quantization of $b$-symplectic manifolds]{On geometric quantization of $b$-symplectic manifolds}

\author{Victor W. Guillemin}
\author{Eva Miranda}
\thanks{{ E. Miranda is partially supported  by the grants reference number MTM2015-69135-P (MINECO/FEDER) and reference number 2014SGR634 (AGAUR)}}
\author{Jonathan Weitsman}
\thanks{J. Weitsman was supported in part by NSF grant DMS 12/11819}
\address{Department of Mathematics, MIT, Cambridge, MA 02139}
\email {vwg@math.mit.edu}
\address{{Department of Mathematics}, Universitat Polit\`{e}cnica de Catalunya and BGSMath, Barcelona, Spain}
\email{eva.miranda@upc.edu}
\address{Department of Mathematics, Northeastern University, Boston, MA 02115}
\email{j.weitsman@neu.edu}
\thanks{\today}

\begin{abstract} We study a notion of pre-quantization for $b$-symplectic manifolds.  We use it to construct a formal geometric quantization of $b$-symplectic manifolds equipped with Hamiltonian torus actions with nonzero modular weight.  We show that these quantizations are finite dimensional $T$-modules.

\end{abstract}
\maketitle

\section{Introduction}

Let $(M,\omega)$ be an integral symplectic manifold, and let $(\cL,\nabla)$ be
a line bundle $\cL$ with connection $\nabla$ of curvature $\omega$.  The
quadruple $(M,\omega,\cL,\nabla)$ is called a {\em prequantization} of $(M,\omega)$, which morally should give rise to a geometric quantization $Q(M)$ of $M$.  A complication arises in that all known constructions of $Q(M)$ require additional data, a polarization
of $M$; such a polarization may be real, a foliation of $M$ by Lagrangian
subvarieties, or else complex, a complex or almost complex structure on $M$
compatible with $\omega.$  It is generally believed, and in many cases verified, that the
quantization $Q(M)$ should be independent of the polarization.  However there
is no theorem guaranteeing that this should be the case.

Work of Kontsevich \cite{Ko} extending deformation quantization to Poisson manifolds
raises the issue as to whether any of the constructions above
has any relevance in the Poisson setting.  If $(M,\pi)$ is a Poisson manifold,
it is not clear what the analog of $(\cL,\nabla)$ should be, let alone
what one would mean by a polarization.  The purpose of this paper is to try
to begin developing some examples, guided by symplectic geometry, where a
sensible theory of geometric quantization of Poisson manifolds can be
proposed.  Hopefully the repertoire of examples may be a guide to a theory
of geometric quantization of Poisson manifolds.

To do this we focus on a special class of Poisson manifolds that have two
helpful properties.  First, we require that the Poisson structure be
symplectic on the complement of a real hypersurface $Z \subset M$ and
have a simple zero on $Z.$  Such $b$-symplectic manifolds have been the
subject of intensive study \cite{GMP, GMPS} and are by now well
understood geometrically.  And second, we require that the manifold have
a Hamiltonian action of a torus with a certain nondegeneracy condition
(nonzero modular weight; see Theorem \ref{gmpsthm} below for the precise definition).  The presence
of these two conditions allows us to bring tools from symplectic
geometry to bear on the problem.   One concept which, as far as we know, has not been investigated at all in the $b$-symplectic setting and which will play an essential ingredient in  describing how to quantize these manifolds, is the concept of "integrality" for the $b$-sympletic form $\omega$ (or, alternatively of "pre-quantizability" for the pair, $(M,\omega)$); and one of the main goals of this paper will be to provide an appropriate definition of this concept and explore some of its consequences.
We then show that a natural functoriality
condition for quantization (``formal geometric quantization'') determines
what the quantization of the manifold should be.

Formal geometric quantization was studied in \cite{Wei} in the
context of the quantization of Hamiltonian $T$-spaces with proper moment map.
We will see that where $M$ is a compact $b$-symplectic manifold, with a Hamiltonian
torus action of nonzero modular weight,
the manifold $M - Z$ is such a space, and that an analog of formal
geometric quantization for $b$-symplectic manifolds yields essentially
the quantization of $M - Z.$  However, in the $b$-symplectic case, there
is a surprise; unlike in the case of noncompact manifolds with proper moment map, where the quantization is always infinite dimensional, though with
finite multiplicities, in the case of a $b$-symplectic manifold, the
quantization is always a {\em finite dimensional} virtual $T$-module.  This
raises the question of whether it is the index of a Fredholm operator.

\section{{$b$}-symplectic manifolds}

Let $M$ be a compact, connected, oriented $n$-dimensional manifold, $Z\subseteq M$ a closed hypersurface and $f:\, M\rightarrow\mathbb{R}, f{|_{Z}}=0$, a defining function for $Z$. We recall (see \cite{GMP}) that a $b$-symplectic form on $M$ is a $2$-form of the form

\begin{equation}
\omega=\frac{df}{f}\wedge\mu+\gamma
	\label{eq:2.1}
\end{equation}

\noindent
with $\mu\,\in\,\Omega^{1}(M)$ and $\gamma\in\Omega^2(M),$ which is symplectic in the usual sense on $M-Z,$ and is symplectic at $p \in Z$ as an element of $\wedge^2(^bT^*_p)$ where $^bT^*_p$ is the span of $T^*_p Z$ and the ``$b$-form'' $\left(\dfrac{df}{f}\right)_p$.

Some properties of the form (\ref{eq:2.1}) which we will need below are:

\begin{enumerate}
	\item Let $\iota: Z\rightarrow M$ be the inclusion map. Then $\iota^*\mu=:\mu_Z$ is an intrinsically defined one-form on $Z$
	\item $f$ is not intrinsically defined but replacing $f$ by $f=hg$ with $h>0$ on $Z$, $\dfrac{df}{f}=\dfrac{dg}{g}+\dfrac{dh}{h}$ so $\dfrac{df}{f}$ {\em is} intrinsically defined $\mod\Omega^1(M)$. Moreover
	
\begin{equation}
\omega=\frac{dg}{g} \;\wedge~\;\mu+\gamma'
	\label{eq:2.2}
\end{equation}

where

\begin{equation}
\gamma'=\gamma+\frac{dh}{h} \;\wedge~\;\mu
	\label{eq:2.3}
\end{equation}
\item Since $d\omega=0=-\dfrac{df}{f}\;\wedge~\;d\mu+d\gamma$ the forms $\iota_Z^*\mu=\mu_Z$ and $\iota_Z^{{*}}\gamma=\gamma_C$ are closed.
\item For $\omega_p, p\,\in\, Z$, to be symplectic in the sense described above, $\mu_Z$ has to be nonvanishing on $Z$ and hence, by item $3$, defines a foliation of $Z$. Moreover {it} also requires that if $L$ is a leaf of this foliation $\iota^*_L\,\gamma$ is a symplectic form on $L$. In addition, by (\ref{eq:2.2}) $\iota^*_L\gamma'={\iota^*_L\gamma}$  so this symplectic structure on $L$ is intrinsically defined.
\end{enumerate}

Turning next to ``pre-quantization'' we note that, since $\mu_Z$ is intrinsically defined, so is its cohomology class, $[\mu_Z]$ and by (\ref{eq:2.3}) the cohomology class $[\gamma]$ is intrinsically defined as well. Moreover the Melrose-Mazzeo isomorphism

\begin{equation}
	^bH^2(M,\mathbb{R})\rightarrow H^2(M,\mathbb{R})\oplus H^1(Z,\mathbb{R})
	\label{eq:2.4}
\end{equation}

\noindent
maps $[\omega]$ onto $[\gamma]\oplus[\mu_Z]$, hence a natural definition of ``integrality'' for $\omega$, i.e. of ``$[\omega] \in ^b\!H^2\,(M,\mathbb{Z})$'' is to require that $[\mu_Z]$ be in $H^1(Z,\mathbb{Z})$ and $[\gamma]$ be in $H^2(M,\mathbb{Z})$. We will list a few consequences of {this} assumption.

\begin{enumerate}
	\item The integrality of $\gamma$ implies that there exists a circle bundle,
	\begin{equation}
	\pi: V\rightarrow M
	\label{eq:2.5}
\end{equation}

and a one form $\alpha$ on $V$ such that

	\begin{equation}
	d\alpha=\pi^*\gamma,
	\label{eq:2.6}
\end{equation}

	\begin{equation}
	\iota (X)\alpha=1
	\label{eq:2.7}
\end{equation}

\qquad and
	\begin{equation}
	{\mathcal{L}}_X\alpha=0
	\label{eq:2.8}
\end{equation}

\noindent where $X$ is the generator of the circle action on $V$.
\item The integrality of {$\mu_{Z}$} implies that there exists a map, \textbabygamma$: Z\rightarrow S^1$, with the property
	\begin{equation}
	\mu_Z=\text{\textbabygamma}^*d\theta
	\label{eq:2.9}
\end{equation}
\end{enumerate}

Therefore, in particular, the foliation of $Z$ that we described above is that defined by the level sets of \textbabygamma, and hence since $Z$ is compact, the leaves of this foliation are compact as well. Moreover if we let {$v$} be the vector field on $Z$ defined by

\begin{equation}
	\iota_v\mu_Z=1 \text{ and } \iota_v\gamma_\iota=0
	\label{eq:2.10}
\end{equation}

then $v$ and $\dfrac{\partial}{\partial\theta}$ are \textbabygamma-related. Therefore if we let $\phi: Z\rightarrow Z$ be the map $\exp 2\pi v$ and let $L$ be a leaf of the foliation defined by \textbabygamma, $Z$ can be identified with the mapping torus

\begin{equation}
	L\times [0, 2\pi]/\sim
	\label{eq:2.11}
\end{equation}

where ``$\sim$'' is the identification

\begin{equation}
	{(p,0) \sim (\phi(p), 2\pi)}
	\label{eq:2.12}
\end{equation}

\section{Group actions}

As in \S $2$ let $f: (M, Z)\rightarrow (\mathbb{R},0)$ be a defining function for $Z$ and let $Z_i, i=1, 2,\ldots, k$, be the connected components of $Z$. We will denote by $^bC^\infty(M)$ the space of functions which are $C^\infty$ on $M-Z$ and near each $Z_i$ can be written as a sum,

\begin{equation}
	c_i\log|f| +g
	\label{eq:3.1}
\end{equation}

\noindent
with $c_i\in\mathbb{R}$ and $g\in C^\infty(M)$.
\def\ft{\mathfrak t}
Now let $T$ be a torus and $T\times M\rightarrow M$ an action of $T$ on $M$.\footnote{The material in this section is taken more or less verbatim from \cite{GMPS}.}
 We will say that this action is {\em Hamiltonian} if the elements, $X\,\in\,\ft$ of the Lie algebra of $T$ satisfy

\begin{equation}
	\iota (X_M)\,\omega=d\phi, \phi~\in~^bC(M),
	\label{eq:3.2}
\end{equation}

\noindent
in other words:

\begin{equation}
	\iota (X_M)\omega=c_i(X)d(\log{|f|}) +d g
	\label{eq:3.2$'$}
\end{equation}

\noindent
in a tube neighborhood of $Z_i$ for
$g\in C^\infty(M)$.

The map

\begin{equation}
	v_i:~ X~\in~ \ft\rightarrow c_i(X)
	\label{eq:3.3}
\end{equation}

\noindent
is called the modular weight of $Z_i$ and depends on $i$; however, one can show (\cite{GMPS}, \S $2.3$)

\vspace{10pt}
\noindent
\begin{RefTheorem}[{\cite{GMPS}}]\label{gmpsthm} The $v_i$'s are either zero for all $i$ or non-zero for all $i$.\end{RefTheorem}

In this paper we will assume that the latter is the case, in which case, for fixed $i$ we can choose {$X_i\,\in\,\mathfrak t$} such that $c_i(X_i)=1$. Hence, by (\ref{eq:2.11}) $\exp 2\pi X_i$ maps the leaves, $L_i$, of the null foliation of $Z_i$ onto themselves, and thus $\exp 2\pi X_i=\exp 2\pi Y_i$ where $Y_i\,\in\,t$ is tangent to the leaves of this foliation. Thus, replacing $X_i$ by $X_i-Y_i$ we can assume that $\exp 2\pi X_i$ is the identity map on $Z$. In other words the map

$$S^1\times L\rightarrow Z, (\theta, p)\rightarrow (\exp\theta{ X_i}) p$$

\noindent
is a diffeomorphism, and hence the mapping tori (\ref{eq:2.11}) are all products: $L\times S^1$. Moreover if we split $T$ into a product

$$T=T_i\times S^1$$

\noindent where $T_i$ is the subgroup of $T$ fixing the leaves of the null-foliation of $Z_i$, the action of $T$ on $Z_i$ is just the product of the canonical action of $S^1$ on $S^1$ and of $T_i$ on $L$.

\section{Formal geometric quantization}

\subsection{Compact symplectic manifolds}  Let $(M,\omega)$ be a compact
symplectic manifold and let $(\cL, \nabla)$ be a line bundle with connection
of curvature $\omega.$  Choose an almost complex structure $J$ compatible
with the symplectic structure.  Then this almost complex structure gives
$\cL$ the structure of a complex line bundle, and by twisting the spin-$\mathbb{C}$ Dirac operator on $M$ by $\cL$ we obtain an elliptic operator $\bar{\partial}_\cL.$ Since $M$ is compact, $\bar{\partial}_\cL$
is Fredholm, and we define the geometric quantization $Q(M)$ by

$$Q(M) = {\rm ind}(\bar{\partial}_\cL)$$

\noindent as a virtual vector space.

If $M$ is equipped with a Hamiltonian action of a torus $T$, the action lifts
to $\cL$, and one can choose the almost complex structure to be $T$-invariant.
Then the quantization $Q(M)$ is a finite-dimensional virtual $T$-module,
and it satisfies the following principle.

For $\xi \in \ft^*,$ we denote by $M//_\xi T$ the reduced space of $M$ at $\xi.$  Also, for $\alpha$ a weight of $T,$ and $V$ a virtual $T$-module, denote by $V^\alpha$ the submodule of $V$ of weight $\alpha.$

\begin{RefTheorem}[\cite{mein}]\label{qr}
Let $\alpha$ be a weight of $T.$  Then
\begin{equation}\label{qreq}
Q(M)^\alpha = Q(M//_\alpha T).\end{equation}
\end{RefTheorem}

In other words,
\begin{equation}\label{qrcons}
Q(M) = \bigoplus_\alpha Q(M//_\alpha T) \alpha
\end{equation}

\noindent as virtual $T$-modules.

\begin{Remark}\label{singrmk} Both Theorem \ref{qr} and equation (\ref{qrcons}) are strictly
speaking valid only for regular values of the moment map.  In the case where
$\alpha$ is a singular value of the moment map, the singular quotient must be replaced
by a slightly different construction using a shift of $\alpha$.  For details, we refer the interested reader to \cite{mein}. A similar caution applies in
the case of noncompact Hamiltonian $T$-spaces and of $b$-symplectic manifolds
below. \end{Remark}
\begin{Remark} If $(M,\omega)$ is a compact, integral symplectic manifold, one can always find a line bundle $\cL$ with connection $\nabla$ of curvature $\omega,$ and the quantization
$Q(M)$ is independent of this choice.  We therefore suppress the line bundle and connection and simply write $Q(M)$ for the quantization.\end{Remark}

\subsection{Noncompact Hamiltonian $T$-spaces}  If we now consider the case
where $M$ is not compact, the analysis above cannot be carried out, since
the operator $\bar{\partial}_\cL$ is elliptic, but no longer Fredholm.  Instead,
in \cite{Wei} (see also \cite{p}), equation (\ref{qreq}) is used to {\em define} the quantization of such
Hamiltonian $T$-spaces, where the moment map is proper, so that the
reduced spaces are compact and the right hand side of equation (\ref{qreq})
makes sense.\footnote{It is also possible in this case to use index theory to define the
quantization; see \cite{b,p}}

\begin{Definition}[\cite{Wei}]\label{fgq} Let $M$ be a Hamiltonian $T$-space with integral symplectic
form.  Suppose the moment map for the $T$-action is proper.
Let $V$ be an infinite-dimensional virtual $T$-module with finite
multipliticies.  We say

$$V= Q(M)$$

\noindent if for any compact Hamiltonian $T$-space $N$ with integral symplectic form, we have

\begin{equation}\label{qreqn1}
(V\otimes Q(N))^T = Q((M \times N)//_0T).\end{equation}

In other words, as in (\ref{qrcons}),

$$Q(M) = \bigoplus_\alpha Q(M//_\alpha T) \alpha,$$

\noindent where the sum is taken over all weights $\alpha$ of $T.$\footnote{Again, care must be taken
about singular values}
\end{Definition}

Note that the fact that the moment map is proper implies that the reduced space $(M \times N)//_0T$ is compact for
any compact Hamiltonian $T$-space $N$, so that the right hand side of equation (\ref{qreqn1}) is well-defined.

In other words, we have used Theorem \ref{qr}  to give us enough functoriality
to force a definition of the quantization in this case, despite the fact that the elliptic operator
$\bar{\partial}_\cL$ is not Fredholm.

\subsection{$b$-symplectic manifolds}

Suppose now that $M$ is a compact $b-$symplectic manifold, with integral
$b-$symplectic form as above.  Suppose that it is equipped with a Hamiltonian
action of a torus $T$ with {\em nonzero modular weight}.  Then, in
analogy with Definition \ref{fgq}, we define

\begin{Definition}
Let $V$ be a virtual $T$-module with finite
multipliticies.  We say

$$V= Q(M)$$

\noindent if for any compact Hamiltonian $T$-space $N$ with integral symplectic form, we have

\begin{equation}\label{qreqn2}
(V\otimes Q(N))^T = Q((M \times N)//_0T).\end{equation}
\end{Definition}

In other words,

$$Q(M) = \bigoplus_\alpha Q(M//_\alpha T) \alpha,$$

\noindent where the sum is taken over all weights $\alpha$ of $T$.\footnote{Again, adjusting for singular values as described in Remark \ref{singrmk}.}
In this $b$-symplectic case the condition that the modular weight be nonzero guarantees that the reduced space $(M \times N)//_0T$
is a compact and symplectic (and in the generic case, a manifold) for any compact Hamiltonian $T$-space $N$; so that
as in the case of noncompact Hamiltonian $T$-spaces, the right hand side of equation (\ref{qreqn2}) is well-defined.

Another way to say this is to note that

$$Q(M) = Q(M - Z)$$

\noindent where $Q(M-Z)$ is the quantization of the noncompact Hamiltonian $T$-space $M-Z$. The fact that the modular weights on $M$ are nonzero insures that the moment map on $M-Z$ is proper.

The main result of this paper is that $Q(M)$ is a {\em finite} virtual $T-$ module.
To see this, we must return to the geometry of the manifold $M$ in the
vicinity of the hypersurface $Z.$

\section{Symmetry properties}

We have shown above that if the modular weight $v_i$ of $Z_i$ is non-zero then, in the vicinity of $Z_i$, $M$ is just a product

\begin{equation}
	Z_i\times (-\epsilon, \epsilon)
	\label{eq:4.1}
\end{equation}

\noindent and {$Z_i=S^1\times L$}.

\noindent
We will show below that this can be slightly strengthened (see also \cite{GMPS0}): the $b$-symplectic form on $Z\times (-\epsilon, \epsilon)$ can be taken to be the two-form

\begin{equation}
	-d\theta\wedge\dfrac{dt}{t} +\gamma_L
	\label{eq:4.2}
\end{equation}

\noindent
where $\gamma_L$ is the symplectic form on $L$ and $-d\theta \wedge \dfrac{dt}{t}$ the standard $b$-symplectic form on $S^1\times (-\epsilon, \epsilon)$.

To see this, we note that, under the hypotheses above, we can assume that the symplectic form (\ref{eq:2.1}) has the form

\begin{equation}
	d\theta \wedge \dfrac{dt}{t}+\gamma_L+d\theta \wedge \beta
	\label{eq:4.5}
\end{equation}

Moreover if $\iota \left(\dfrac{\partial}{\partial\theta}\right)\beta=h$ we can replace $\beta$ by $\beta-h~d\theta$ in the expression above and arrange that $\iota \left(\dfrac{\partial}{\partial\theta}\right)\beta=0$. Hence since the action of $S^1$ on $M$ is Hamiltonian

$$\iota{\left(\dfrac{\partial}{\partial\theta}\right)}\omega=d(\log |t|+{\rho})$$

\noindent for some $\rho~\in~ C^\infty(M)$ and hence

\begin{equation}
	\beta=d{\rho}
	\label{eq:4.6}
\end{equation}

Consider now the one parameter family of forms

\begin{equation}
	d\theta \wedge \dfrac{dt}{t}+\gamma_L-s d({\rho} d\theta)
	\label{eq:4.7}
\end{equation}

\noindent
for $0\leq s\leq 1$. For $s=1$ this form is $\omega$ and for $s=0$ the form (\ref{eq:4.2}). Moreover for $\epsilon$ small and $-\epsilon<t<\epsilon$ the first summand of {(\ref{eq:4.7})} is much larger than the third so the form {(\ref{eq:4.7})} is $b$-symplectic and for all $s$

$$[\omega_s]=[\omega_0]$$

\noindent
so we can apply $b$-Moser theorem (see \cite{GMP}) to conclude that $\omega_0$ and $\omega_1$ are equivariantly symplectomorphic.

\vspace{10pt}
{Finally note that} the $2$-form, $\gamma_L$, depends in principle on $t$.

However the inclusion map

$$\iota: L\rightarrow L\times(-\epsilon, \epsilon),~p\rightarrow (p,0)$$

and the projection map

$$\pi: L\times (-\epsilon, \epsilon)\rightarrow L, (p, e)\rightarrow p$$

\noindent
induce isomorphisms on cohomology; hence {$[\mu]=[\pi^* \iota^*\mu]$}. Therefore, by Moser, we can assume that {$\mu=\pi^* \iota^* \mu$} i.e. that $\mu$ is just a symplectic $2$-form on $L$ itself.

\section{Formal quantization of $b$-symplectic manifolds}

We now prove the main result of this paper.

\begin{Theorem}\label{finite}

Let $M$ be an integral $b$-symplectic manifold equipped with a Hamiltonian $T$-action with nonzero modular weight.
Then the formal geometric quantization $Q(M)$ is a finite dimensional $T$-module.\end{Theorem}

\begin{Proof}
We will show that if we take for the quantization of $N=Z_u\times (-\epsilon, \epsilon)$ the sum:

\begin{equation}
	\oplus~Q(N//_\alpha T) \alpha,~\alpha~\in~\mathbb{Z}~(T)
	\label{eq:4.4}
\end{equation}

\noindent
then the virtual vector spaces $Q(N//_\alpha T) $, are all zero and hence so is this sum.

To define $Q(N//_\alpha T)$ one has to define orientations on the $N//_\alpha T$ and to do this consistently one has to define an orientation on $N$. The natural way to do so would be to assign to each connected component of $N$-$Z$ the orientation defined by the symplectic form, $\omega$; however, because of the factor $\dfrac{df}{f}$, in the formula (\ref{eq:2.1}) the symplectic orientations on adjacent components of the space $N$-$Z$ don't agree; and, in particular, if $N = Z_i\times (-\epsilon, \epsilon)$ and $\mathbb{R}\times \mathfrak t_i$ is the Lie algebra of $S^1\times T_i$ the moment map

{$$\phi : N\rightarrow \mathbb{R}\times  \ft^*_i$$}

\noindent
associated with the action of $S^1\times T_i$ on $N$ is the map

$$(\theta, p, t)~\epsilon~S^1\times L\times(-\epsilon, \epsilon)\rightarrow {(\log|t|,{~\phi_i (p))}}$$

\noindent
where $\phi_i: L\rightarrow \ft^*_i$ is the moment map associated with the $T_i$ action on $L$. Thus for a weight, $(c, \alpha_i)$ of $S^1 \times T_i$, the reduced space

$$\phi^{-1}(c, \alpha_i)/S^1\times T_i$$

\noindent
consists of two copies of the reduced space

$$\phi^{-1}_i(\alpha_i)/T_i$$

\noindent
with opposite orientations, so the quantization of this space is a virtual vector space $V_+\oplus V_-$ with $V_-=-V_+$. \end{Proof}

We end the paper with a conjecture.

\begin{Conjecture}  There exists a natural Fredholm operator on $M$ whose index gives $Q(M)$.\end{Conjecture}

\vspace{5pt}
\noindent

\end{document}